\documentstyle{amsppt}
\magnification=1000
\NoBlackBoxes

\def\R{\Bbb R}

\def\back{\backslash}
\def\fp{\flushpar}
\def\g{\gamma}
\def\lra{\longrightarrow}
\def\mt{\mu_\theta}
\def\mgt{\mu_{\beta,\theta}}
\def\pl{\pi_L}
\def\pr{\pi_R}
\def\rson{\R^n\times SO(n)}
\def\SN{\Bbb S^{n-1}}

\def\txn{T^*(\R^n\times \R^k)}

\def\tyn{T^*\R^n}
\def\G{\Gamma}

\loadbold\topmatter
\title
$L^p-L^{p'}$ estimates for overdetermined Radon transforms
\endtitle
\author
Luca Brandolini, Allan Greenleaf and Giancarlo Travaglini
\endauthor
\thanks 
The second author was partially supported by a grant from the  National
Science Foundation.
 \endthanks

\address
Dipartimento di Ingegneria Gestionale e dell'Informazione, Universit\`a degli Studi
di Bergamo, V.le G Marconi 5, 24044 Dalmine, ITALY  
\endaddress
\email brandolini\@unibg.it\endemail

\address
Department of Mathematics, University of Rochester, Rochester, NY 14627
\endaddress
\email allan\@math.rochester.edu \endemail
\address
Dipartimento di Matematica e Applicazioni, Universit\`a di Milano-Bicocca, Via
Bicocca degli Arcimboldi 8, 20126 Milano, ITALY  
\endaddress
\email giancarlo.travaglini\@unimib.it\endemail

\abstract
We prove several variations on the results in Ricci and Travaglini[RT] 
concerning
$L^p-L^{p'}$ bounds for convolution with all rotations of a measure supported by
a fixed convex curve in $\R^2$. Estimates are obtained for averages over
higher-dimensional convex (nonsmooth) hypersurfaces, smooth $k$-dimensional
surfaces, and nontranslation-invariant families of surfaces. We compare the
approach of [RT],  based on average decay of the Fourier transform,
with an approach based on $L^2$ boundedness of Fourier integral operators, and show
that essentially the same geometric condition arises in proofs using 
different techniques.
\endabstract

\endtopmatter

\head{\bf \S1. Introduction}\endhead

Our starting point is the following result from [RT]:
\bigskip

\proclaim{Theorem 1} Let $\G\subset\R^2$ be a compact, convex curve with
arc length measure $\mu$. Let $\mt$ denote the rotation of $\mu$ by
$\theta\in \Bbb S^1$. Then
$$\Bigl(\int_{\Bbb S^1}\int_{\R^2} |f*\mt (x)|^3 dx d\theta \Bigr)^\frac13
\lesssim\Bigl(\int_{\R^2} |f(y)|^\frac32 dy \Bigr)^\frac23.\tag1.1$$
\endproclaim
\flushpar Thus, the
$L^{\frac32}(\R^2)\lra L^3(\R^2)$ estimate that holds for curves in the plane with
nonzero curvature ([L],[Str]) generalizes to arbitrary (i.e., not necessarily
smooth) convex curves when averaged over all rotations. The goal here is to extend
this in several ways: to averages over
$k$-dimensional surfaces in $\R^n$; to more general transformations than
rotations; and to nontranslation-invariant averaging operators. In doing
so, we will primarily use two techniques: estimates for average decay of
the Fourier transform of surface measure (as in [RT]), and $L^2$
regularity of nondegenerate Fourier integral operators. Although these
methods appear to be  different, the geometric assumptions needed to use
them are essentially the same. 

To start with, one can extend Thm. 1 to rotations of  hypersurfaces in
higher dimensions with the same convexity assumption. For $\theta\in SO(n)$ and
$\mu$ a measure on
$\R^n$, let $\mt$ be defined by $<f,\mt>=<f(\theta^{-1}\cdot),\mu>$, so
that $\hat\mu_\theta(\xi)=\hat\mu(\theta\xi)$.

\bigskip
\proclaim{Theorem 2} Let $S\subset\R^n$ be a compact, convex hypersurface
with induced measure $\mu$. Then,

$$\Bigl(\int_{SO(n)}\int_{\R^n} |f*\mt (x)|^{n+1} dx
d\theta\Bigr)^\frac{1}{n+1}
\lesssim \Bigl(\int_{\R^n} |f(y)|^\frac{n+1}{n} dy
\Bigr)^\frac{n}{n+1}.\tag1.2$$

\endproclaim
\bigskip

(Here, and throughout, we use $\lesssim$ to denote $\le c\cdot$, with $c$
dependent only on the operator in question.) 
Thm. 1 can also be modified to cover all rotations of a surface in $\R^n$
of arbitrary dimension, under a smoothness assumption. 

\bigskip
\proclaim{Theorem 3} Let $S\subset \R^n$  be a smooth $k$-dimensional
surface, $1\le k\le n-1, n\ge 2$, and $\mu$ a smooth, compactly supported
multiple of induced surface measure on $S$. Then,

$$\Bigl(\int_{SO(n)}\int_{\R^n} |f*\mt (x)|^{\frac{2n-k}{n-k}} dx
d\theta\Bigr)^\frac{n-k}{2n-k}
\lesssim\Bigl(\int_{\R^n} |f(y)|^\frac{2n-k}{n} dy
\Bigr)^\frac{n}{2n-k}.\tag1.3$$

\endproclaim
\bigskip

A crucial ingredient in the proof of Thm. 1 was the $L^2$ average decay of the
Fourier transform $\hat\mu$ from [P]. Thms. 2 and 3 follow immediately by
replacing Podkorytov's estimate in the argument of [RT] by the results of [BHI] and
Prop. 1 below, respectively.

To obtain the optimal $L^\frac{2n-k}{n}\lra
L^\frac{2n-k}{n-k}$ boundedness, we do not actually need to use all rotations of
the surface, or even linear transformations for that matter, nor does the operator
need to be translation-invariant. To start with, we keep the
translation-invariance, but allow nonlinear transformations to act on the
surface.  Let
$T_{s}:\R^{n}\rightarrow\R^{n}$
be a smooth family of transformation of $\R^{n}$ parametrized by $s\in
K\subset\subset\R^{m},\quad m\ge n-1$, let $\gamma
:A\subset\R^{n-1}\rightarrow\R^{n}$ be a $C^{2}$ parametrized convex
hypersurface in $\R^{n}$ (if $n=2$ we merely need to assume that $\g$ is convex)
and set
$\gamma_{s}\left(  t\right) =T_{s}\left(  \gamma\left(  t\right) 
\right)  $. We are interested in the operator
$$
Rf\left(  x,s\right)  =\int f\left(  x+\gamma_{s}\left(  t\right)
\right) \chi(t) dt=(f\ast\mu_{s})\left(  x\right)
$$
where $\mu_{s}$ is the measure defined by
$$
\int f\left(  x\right)  d\mu_{s}=\int f\left(  \gamma_{s}\left(
t\right)  \right)  \chi(t) dt,
$$
with $\chi\in C^\infty_0(\R^{n-1})$ is a fixed cutoff function.
Denote by $J_{T_{s}}\left(  x\right)  $ and $J_{\frac{\partial
T_{s}}{\partial\theta_{k}}}\left(  x\right)  $ the Jacobian matrices at
$x$ of
the maps $T_{s}$ and $\frac{\partial T_{s}}{\partial s_{k}}$,
respectively. We have
\bigskip
\proclaim{Theorem 4} Let $\g$ be a convex curve if $n=2$ and a $C^2$ convex
hypersurface if $n\ge 3$. Assume that for every unit vector $\Xi\in\SN$ and for
every
$s\in K$ the matrix
$$
C=\left [
\matrix
\Xi^{t}J_{T_{s}}\left(  x\right) \cr
\Xi^{t}J_{\frac{\partial T_{s}}{\partial s_{1}}}\left(  x\right)
\cr
\vdots\cr
\Xi^{t}J_{\frac{\partial T_{s}}{\partial s_{n-1}}}\left(
x\right)
\endmatrix
\right ]\tag1.4
$$
has rank $n$. Then
$$
\left\{  \int_K\int_{\R^n}\left\vert Rf\left(  x,s\right)  \right\vert
^{n+1}dxd s\right\}  ^{1/\left(  n+1\right)  }\lesssim\left\{
\int_{\R^n}\left\vert f\left(  x\right)  \right\vert ^{\left(  n+1\right)
/n}dx\right\}  ^{n/\left(  n+1\right)  }\;.\tag1.5
$$
\endproclaim
\bigskip

To consider nontranslation-invariant operators, we now take
$$\gamma:\R^n\times \R^m\times \R^k\lra \R^n$$
to be a $C^\infty$ map, with $D_t\gamma(x,s,t)$ injective, so that each
$\gamma_{x,s}:=\{\gamma(x,s,t): t\in\R^k\}$ is a smooth immersed
$k$-surface in $\R^n$. Define the (overdetermined) generalized Radon
transform $R:C^\infty(\R^n)\lra C^\infty(\R^n\times \R^m)$ by
$$Rf(x,s)=\int_{\R^k} f(\gamma(x,s,t))\chi(t) dt, \quad \chi\in
C_0^\infty(\R^k).\tag1.6$$

Then (see [Gu,GuSt]) the Schwartz kernel of $R$ is a smooth density
$\delta_Z$ supported on the incidence relation
$$Z=\{(x,s,y): y\in\g_{x,s}\}\subset\R^{n+m}\times\R^n,$$
which is codimension $n-k$ in $\R^{n+m}\times\R^n$. If
$$Z=\{F_1(x,s,y)=\dots =F_{n-k}(x,s,y)=0\}$$ 
locally, with $\{dF_1,\dots,dF_{n-k}\}$ linearly independent, then $\delta_Z$ has
the oscillatory representation
$$\delta_Z=\int_{\R^{n-k}} e^{i\sum_{j=1}^{n-k} F_j(x,s,y)\theta_j}
a(x,s,y;\theta) d\theta\tag1.7$$
in the sense of [H1], with $a(x,s,y;\theta)$ a symbol of order 0 (essentially
$\equiv 1$ in $\theta$). In general, $\delta_Z$ is a locally finite sum of such
expressions.  Thus,
$\delta_Z$ is a Fourier integral distribution on
$\R^{n+m}\times\R^n$ associated to the conormal bundle of $Z$,
$$N^*Z=\Bigl\{(x,s,\xi,\sigma;y,\eta)\in T^*\R^{n+m}\times T^*\R^n:
(x,s,y)\in Z, (\xi,\sigma,\eta)\perp TZ\Bigr\}$$
and hence
$R$ is a Fourier integral operator,
$$R\in I^r(\R^{n+k},\R^n;C),$$
where 
$$C=N^*Z'=\Bigl\{(x,s,\xi,\sigma;y,\eta): (x,s,\xi,\sigma;y,-\eta)\in
N^*Z\Bigr\}$$
is a canonical relation, i.e., a lagrangian submanifold for the
difference symplectic form $\omega_{T^*\R^{n+m}}-\omega_{T^*\R^n}$ on
$T^*\R^{n+m}\times T^*\R^n$, and  the order
$r$ is calculated by 
$$\eqalign{r=&(\hbox{ order of } a) +\frac{\hbox{
number of phase variables}}{2}-\frac{\hbox{ number of
spatial variables}}4 \cr
=& 0 +\frac{n-k}2 -\frac{n+m+n}4=-\frac{k}2 -\frac{m}4.}$$

$L^2$ estimates for Fourier integral operators associated with a
canonical relation $C\subset T^*X\times T^*Y$ depend on the structure of
the projections $\pr:C\lra T^*Y$ and $\pl:C\lra T^*X$. The optimal $L^2$
estimates for  an operator $F\in I^r(X,Y;C)$  hold under the assumption that $\pr$
is a submersion (which guarantees that $\pl$ is an immersion), together with the
mild requirement that the spatial projections $\pi_X:C\lra X$ and $\pi_Y:C\lra Y$
are submersions [H1,H2]; such canonical relations
$C$ are called {\it nondegenerate}. Substituting $L^2$ estimates for such
operators in place of the average decay estimates for the Fourier transform of
surface-carried measures, we can show

\bigskip
\proclaim{Theorem 5} Let $R$ be a generalized Radon transform as in (1.6)
such that the associated canonical relation $C$ is nondegenerate. Then,
$$Rf:L^{\frac{2n-k}{n}}_{comp}(\R^n)\lra L^\frac{2n-k}{n-k}_{loc}
(\R^{n+m}).\tag1.8$$

\endproclaim
\bigskip
Letting $\g_0:\R^k_t\lra \R^n$ be a local parametrization of a smooth $k$-surface
$S\subset\R^n$, $m=\frac{n(n-1)}2$, $\theta:\R^m_s\lra SO(n)$ a coordinate
chart and
$\g(x,s,t)=x-\theta(s)(\g_0(t))$, we see that Thm. 5 extends Thm. 3. It is also
possible to use proper subgroups of $SO(n)$ and obtain the same estimates. These
and other particular cases of Thm. 5 will be discussed in
\S5 below.

\bigskip

All of these results involve estimates on the line of duality. Via interpolation
with the $L^1-L^1$ and $L^\infty-L^\infty$ bounds, we find that the type sets of
the operators contain certain closed triangles, symmetric about the line of
duality. For general hypersurfaces, this is sharp, as the example of the
unit sphere shows, with rotation not producing any additional $L^p$
improvement. The emphasis here is on the extension of these estimates to low
regularity and variable coefficient settings. For higher codimension
surfaces, the results here fail to be sharp. For example, Drury [D] (see
also Christ [C1]) has shown that the $X$-ray transform on $\R^n$ maps
$L^{\frac{n+1}2}(\R^n)$ to $L^{n+1}(M_{1,n})$, where $M_{1,n}$ is the
Grassmannian of affine lines in $\R^n$, and this then implies an improvement of
Thm. 3 for $S$ a line segment. Also, these
results have a somewhat different character then those of, for example,
[O],[PhS],[GSW],[S],[C2] or [TaW], where the specific geometry of the curve or
family of curves determines a more complicated type set in the absence of
rotations.
\bigskip

Finally, mixed norm estimates are possible for certain
model surfaces in $\R^n$, just as in [RT] for model curves in $\R^2$. Writing
$x=(x',x_n)$, consider the hypersurface

$$
S_{\beta}=\left\{  \left(  x^{\prime},x_{n}\right)  \in\R^{n}
:x_{n}=\left\vert x^{\prime}\right\vert ^{\beta}\right\}
$$
where $\beta>2$. Let $d\mu_\beta$ be the  induced  measure
on $S_{\beta}$, multiplied by a $C^{\infty}$ function on
$\R^{n}$ with compact support, identically $1$ in a neighborhood of
the origin, and $\mgt$ its rotation by $\theta\in SO(n)$. We have
\bigskip
\proclaim{Theorem 6} Let $R_\beta f(x,\theta)=f\ast\mgt(x)$.
Then,
$$
\left\Vert R_\beta f\right\Vert _{L^{\frac{p^{\prime}\left(
n-1\right)  \left(  \beta-1\right)  }{\left(  p^{\prime}-2\right)
\beta-2\left(  n-1\right)  }}\left(  SO\left(  n\right) ;L^{p^{\prime}
}\left(  \R^{n}\right)  \right)  }\lesssim\left\Vert f\right\Vert
_{L^{p}},\quad \frac{n+1}{n}<p<\frac{2\beta+2\left(  n-1\right)  }{\beta
+2\left(  n-1\right)  }\tag{i}$$

$$\left\Vert R_\beta f\right\Vert _{L^{\infty}\left(
SO\left(  n\right) ;L^{p^{\prime}}\left(  \R^{n}\right)  \right)
}\lesssim\left\Vert f\right\Vert _{L^{p}},\quad p>\frac{2\beta+2\left(
n-1\right)  }{\beta+2\left(  n-1\right)}\tag{ii}
$$
\endproclaim

Thms. 2 and 3, which are based on $L^2$ average decay properties of Fourier
transforms of surface-carried measure, are proved in \S2. Thm. 6, which uses $L^p$
average decay properties, is proved in \S3, and Thm. 4, which still concerns
translation-invariant operators and thus can be treated using Fourier transform
estimates, is treated in
\S4.  In \S5, we prove Thm. 5  and  discuss geometric criteria for nondegeneracy
of the canonical relation. 
\bigskip

\head{\bf \S2. Euclidian motions of a fixed surface}\endhead

We begin by considering averages over all translations and rotations of a fixed
$k$-dimensional surface in $\R^n$.
\bigskip

\proclaim{Proposition 1}
Let $\Phi:\Omega\subseteq\R^{k}\rightarrow\R^{n}$ be a
parameterization of a  $C^{\infty}$ $k$-dimensional   surface 
$S\subset\R^{n}$ and let $\mu=\Phi_*(\chi(u) du)$, where $\Phi_*$ denotes 
pushforward and $\chi$ is a suitable cut-off function on $\R^k$, so that
$$\widehat{\mu}\left(  \rho\omega\right)  =\int_{\Omega}e^{-2\pi i\rho\omega\cdot
\Phi\left(  u\right)  }\eta\left(  u\right)  du\tag2.1$$
for $\omega\in\SN$. Then
$$\left\{  \int_{\SN}\left\vert \widehat{\mu}\left(  \rho\omega\right)
\right\vert ^{2}d\omega\right\}  ^{\frac{1}{2}}\lesssim
\rho^{-\frac{k}{2}}.\tag2.2$$

\endproclaim

\fp{\bf Proof}. We can change parameterization and choose coordinates in
$\R^{n}$ so that
$$\Phi\left(  u\right)  =\left(  u,\Psi\left(  u\right)  \right).$$
We can also assume that $\Phi\left(  0\right)  =0$ and $\nabla\Psi\left(  0\right)
=0$. We have
$$\widehat{\mu}\left(  \rho\omega\right)  =\int_{\Omega}e^{-2\pi i\rho\omega
\cdot\left(  u,\Psi\left(  u\right)  \right)  }J\left(  u\right)  du,$$
where $J\left(  u\right) $ is a suitable function that takes into account the
change of parameterization. Let $\psi\left(  \omega\right)  =\psi\left(
\omega_{1},\ldots,\omega_{k}\right) $ be a cut-off function supported on
$\left\vert \omega_{j}\right\vert \leqslant\frac{1}{10}$ for $j=1,\ldots k$.
Then
\medskip
$$\eqalign{
\int_{\SN}\left\vert \widehat{\mu}\right\vert ^{2} & \left(  \rho\omega\right)
\psi\left(  \omega\right)  d\omega\cr  &  =\int_{\Omega}
\int_{\Omega}\int_{\SN}e^{-2\pi i\rho\omega\cdot\left(  u-v,\Psi\left(
u\right)  -\Psi\left(  v\right)  \right)  }\psi\left(  \omega\right)
d\omega\,J\left(  u\right)  J\left(  v\right)  du\,dv\cr
&  =\int_{\Omega}\int_{\Omega}\frac{1}{\left(  1+\rho\left\vert u-v\right\vert
\right)  ^{N}}\,J\left(  u\right)  J\left(  v\right)  du\,dv\cr
&  \lesssim\int_{\Omega}\int_{\R^{k}}\frac{1}{\left(  1+\rho\left\vert
u-v\right\vert \right)  ^{N}}\,du\,J\left(  v\right)  dv\cr
&  =\int_{\Omega}\int_{\R^{k}}\frac{1}{\left(  1+\rho\left\vert
u\right\vert \right)  ^{N}}\,du\,J\left(  v\right)  dv=c \rho^{-k}\int_{\Omega
}\,J\left(  v\right)  dv\cr
&  \lesssim\rho^{-k}\left\vert \Omega\right\vert}
$$
for $N>k$. Here we used the fact that $\omega$ and $\left(  u-v,\Psi\left( 
u\right) -\Psi\left(  v\right)  \right)  $ are almost orthogonal on the support of
$\psi$ and we can evaluate the integral on $\SN$ integrating by parts
$N$ times.

Let $\omega$  be in the support of $1-\psi\left(  \omega\right)  $. Then
$$\widehat{\mu}\left(  \rho\omega\right)  =\int_{\Omega}e^{-2\pi i\rho\omega
\cdot\left(  u,\Psi\left(  u\right)  \right)  }J\left(  u\right)  du.$$
Since $\left\vert \nabla_{u}\left(  \omega\cdot\left(  u,\Psi\left(  u\right)
\right)  \right)  \right\vert >c>0$, integrating by parts $N$ times gives
$$\left\vert \widehat{\mu}\left(  \rho\omega\right)  \right\vert \lesssim
\rho^{-N},$$
finishing the proof of Prop. 1. \qed

\bigskip
\fp{\bf Remark.} Related results for curves in $\R^n$ can be found in [M].
\bigskip

\fp{\bf Proofs of Thms. 2 and  3.} By decomposing $S$ into a finite number of
pieces we can  assume that $S$ is
defined by

$$
\left\{  \left(  x,x^{\prime}\right)  \in\R^{k}\times\R
^{n-k}:x^{\prime}=\Phi\left(  x\right)  \right\}
$$
and that the Jacobian of $\Phi$ has bounded entries. Also observe that the
tangent spaces to $S$ do not contain any line parallel to $\left\{
0\right\}  \times\R^{n-k}$.

Let $i_{z}$ be the distribution defined by

$$
\left\langle i_{z},\varphi\right\rangle =\frac{1}{\Gamma\left(  z\right)
}\int_{0}^{\infty}\varphi\left(  t\right)  t^{z-1}dt,\tag2.3
$$
for test functions $\varphi\in C^\infty_0(\R)$,
and for $\theta\in SO(n)$ let the distribution $\mu_{\theta}^{z}$ be defined by

$$
\widehat{\mu_{\theta}^{z}}\left(  \xi\right)  =\widehat{\mu_{\theta}}\left(
\xi\right)  \widehat{~i_{z}}\left(  \left(  \theta\xi\right)  _{n-k+1}\right)
~\cdots~\widehat{i_{z}}\left(  \left(  \theta\xi\right)  _{n}\right),\tag2.4
$$
where

$$
\Bigl(\left(  \theta\xi\right)_{n-k+1}, \dots,\left(  \theta\xi\right)
_{n}\Bigr)
$$
denotes the last $k$ components of $\theta\xi$. Introducing the analytic
family of operators

$$
T^{z}f(x,\theta)=(f\ast\mu_{\theta}^{z})(x),\quad z\in \Bbb C,\tag2.5
$$
the proof now follows exactly as in [RT]:  using either [BHI] for Thm. 2 or
Prop. 1 for Thm. 3, one shows that

$$
T^{-\frac{k}{2\left(  n-k\right)  }+i\sigma}:L^{2}\left(  \R
^{n}\right)  \rightarrow L^{2}\left(  \R^{n}\times SO\left(  n\right)
\right),\quad\sigma\in\R,\tag2.6
$$
and by (2.4),
$$
T^{1+i\sigma}:L^{1}\left(  \R^{n}\right)  \rightarrow L^{\infty
}\left(  \R^{n}\times SO\left(  n\right)  \right),\quad \sigma
\in\R,\tag2.7
$$
Analytic interpolation then yields that

$$T^0:L^{\frac{2n-k}{n-k}}(\R^n)\lra L^{\frac{2n-k}n}(\rson),$$
which is (1.2) (for $k=n-1$) and (1.3).\qed
\bigskip

\bigskip
\head{\bf \S3. Mixed norm estimates for model surfaces}\endhead

Consider the hypersurface
$$
S_{\beta}=\left\{  \left(  x^{\prime},x_{n}\right)  \in\R^{n}
:x_{n}=\left\vert x^{\prime}\right\vert ^{\beta}\right\}
$$
for values of $\beta>2$. Let $\mu_\beta$ be the measure induced by the Lebesgue
measure on $S_{\beta}$, multiplied by $\widetilde{\psi}\in C^{\infty}_0(\R^n)$
, identically $1$ in a neighborhood of
the origin. We are first interested in the
decay at infinity of the Fourier transform of this measure,

$$
\widehat{\mu_\beta}\left(  \xi\right)  =\int_{\R^{n-1}}e^{-2\pi i\left(
\xi^{\prime}\cdot x^{\prime}+\xi_{n}\left\vert x^{\prime}\right\vert ^{\beta
}\right)  }\psi\left(  x^{\prime}\right)  dx^{\prime}.
$$

\bigskip
\proclaim{Lemma 1}
We have
$$
\left|  \widehat{\mu_\beta}\left(  \xi\right)  \right|  \lesssim \left|
\xi^{^{\prime}}\right|  ^{-\frac{\left(  n-1\right)  \left(  \beta-2\right)
}{2\left(  \beta-1\right)  }}\left|  \xi_{n}\right|  ^{-\frac{n-1}{2\left(
\beta-1\right)  }} \tag3.1
$$
and
$$
\left|  \widehat{\mu_\beta}\left(  \xi\right)  \right|  \lesssim \left|
\xi\right|  ^{-\frac{n-1}{\beta}} \tag3.2
$$

\endproclaim

\fp{\bf Proof.}
To prove (3.1) let $\varphi\left(  x^{\prime}\right)  =\psi\left(
x^{\prime}\right)  -\psi\left(  2x^{\prime}\right)  $. Then,
$$\eqalign{
\widehat{\mu_\beta}\left(  \xi\right)   &  =\sum_{j=1}^{+\infty}\int_{\R
^{n-1}}e^{-2\pi i\left(  \xi^{\prime}\cdot x^{\prime}+\xi_{n}\left|
x^{\prime}\right|  ^{\beta}\right)  }\varphi\left(  2^{j}x^{\prime}\right)
dx^{\prime}\cr
&  =\sum_{j=1}^{+\infty}2^{-\left(  n-1\right)  j}\int_{\R^{n-1}
}e^{-2\pi i\left(  2^{-j}\xi^{\prime}\cdot u+2^{-\beta j}\xi_{n}\left|
u\right|  ^{\beta}\right)  }\varphi\left(  u\right)  du\cr
&  =\sum_{j=1}^{+\infty}2^{-\left(  n-1\right)  j}I\left(  2^{-j}\xi^{\prime
},2^{-\beta j}\xi_{n}\right)}\tag3.3
$$
where
$$
I\left(  \xi\right)  =\int_{\R^{n-1}}e^{-2\pi i\left(  2^{-j}
\xi^{\prime}\cdot u+2^{-\beta j}\xi_{n}\left|  u\right|  ^{\beta}\right)
}\varphi\left(  u\right)  du.
$$

Since $\psi$ is identically $1$ in a neighborhood of the origin, $\varphi$ is
supported away from the origin. Therefore for $I$ we have the estimate
$$
\left|  I\left(  \xi\right)  \right|  \lesssim\left|  \xi\right|
^{-\frac{n-1}{2}}\quad ,
$$
since $S_\beta$ has strictly positive
Gaussian curvature away from the origin. It follows that

$$
\left|  \widehat{\mu_\beta}\left(  \xi\right)  \right|  \lesssim\sum
_{j=1}^{+\infty}2^{-\left(  n-1\right)  j}\left|  \left(  2^{-j}\xi^{\prime
},2^{-\beta j}\xi_{n}\right)  \right|  ^{-\frac{n-1}{2}}.
$$

Let $j_{o}=\left[  \frac{1}{\beta-1}\log_{2}\left(  \frac{\left|  \xi
_{n}\right|  }{\left|  \xi^{\prime}\right|  }\right)  \right]  $ where
$\left[  \cdot\right]  $ denotes the integral part. Observe that when
$j\leqslant j_{o}$ we have

$$
\left|  \left(  2^{-j}\xi^{\prime},2^{-\beta j}\xi_{n}\right)  \right|
\approx2^{-\beta j}\left|  \xi_{n}\right|
$$
while for $j>j_{o}$ we have%

$$
\left|  \left(  2^{-j}\xi^{\prime},2^{-\beta j}\xi_{n}\right)  \right|
\approx2^{-j}\left|  \xi^{\prime}\right|  .
$$
Splitting the above series yields
$$
\eqalign{
\left|  \widehat{\mu_\beta}\left(  \xi\right)  \right|   &  \lesssim\sum
_{j=1}^{j_{o}}2^{-\left(  n-1\right)  j}\left(  2^{-\beta j}\left|  \xi
_{n}\right|  \right)  ^{-\frac{n-1}{2}}+\sum_{j=j_{o}+1}^{+\infty}2^{-\left(
n-1\right)  j}\left(  2^{-j}\left|  \xi^{\prime}\right|  \right)
^{-\frac{n-1}{2}}\cr
&  \lesssim\left|  \xi^{\prime}\right|  ^{-\frac{\left(  \beta-2\right)
\left(  n-1\right)  }{2\left(  \beta-1\right)  }}\left|  \xi_{n}\right|
^{-\frac{n-1}{2\left(  \beta-1\right)  }}}
$$

To prove (3.2) we observe that it is enough to consider the case $\left|
\xi^{\prime}\right|  <c\left|  \xi_{n}\right|  $. From (3.3) we get

$$
\eqalign{
\left|  \widehat{\mu_\beta}\left(  \xi\right)  \right|   &  \lesssim\sum
_{j=1}^{+\infty}2^{-\left(  n-1\right)  j}I\left(  2^{-j}\xi^{\prime
},2^{-\beta j}\xi_{n}\right) \cr
&  \lesssim\sum_{j=1}^{+\infty}2^{-\left(  n-1\right)  j}\left(
1+2^{-\beta j}\left|  \xi_{n}\right|  \right)  ^{-\frac{n-1}{2}}\cr
&  =\sum_{j=1}^{j_{0}}2^{-\left(  n-1\right)  j}\left(  1+2^{-\beta
j}\left|  \xi_{n}\right|  \right)  ^{-\frac{n-1}{2}}+\sum_{j=j_{0}%
+1}^{+\infty}2^{-\left(  n-1\right)  j}\left(  1+2^{-\beta j}\left|  \xi
_{n}\right|  \right)  ^{-\frac{n-1}{2}},}
$$
where $j_{0}=\left[  \frac{1}{\beta}\log\left|  \xi_{n}\right| 
\right]$.
Therefore,
$$
\eqalign{
\left|  \widehat{d\mu}\left(  \xi\right)  \right|   &  \lesssim\sum
_{j=1}^{j_{0}}2^{-\left(  n-1\right)  j}2^{\beta j\frac{n-1}{2}}\left|
\xi_{n}\right|  ^{-\frac{n-1}{2}}+c\sum_{j=j_{0}+1}^{+\infty}2^{-\left(
n-1\right)  j}\cr
&  \lesssim\left|  \xi_{n}\right|  ^{-\frac{n-1}{\beta}}\lesssim\left|
\xi\right|  ^{-\frac{n-1}{\beta}}.}
$$

Lemma 1 allows us to obtain $L^p$ average decay of $\widehat{\mu_\beta}$,
extending a result in [BRT].

\bigskip
\proclaim{Proposition 2}
We have the following estimates:

$$
\Bigl\{  \int_{\SN}\left|  \widehat{\mu_\beta}\left(  \rho\omega\right)
\right|  ^{p}d\omega\Bigr\}  ^\frac{1}{p}\lesssim
\cases\rho^{-\frac{n-1}{2}} & p<\frac{2\left(  \beta-1\right)  }{\beta-2}\cr
\rho^{-\frac{n-1}{2}}\log\left(  \rho\right)  ^{\frac{\left(  \beta-2\right)
\left(  n-1\right)  }{2\left(  \beta-1\right)  }} & p=\frac{2\left(
\beta-1\right)  }{\beta-2}\cr
\rho^{-\left(  n-1\right)  \left(  \frac{1}{p}+\frac{1}{\beta}-\frac
{1}{p\beta }\right)  } & p>\frac{2\left(  \beta-1\right)  }{\beta-2}
\endcases
$$
\endproclaim

\fp{\bf Proof.}
Let now $\xi=(\xi^{\prime},\xi_{n})=(\rho\omega^{\prime}\sin\theta,\rho
\cos\theta),$ with $\omega^{\prime}\in\Bbb S^{n-2}$ and $0\leq\theta\leq\pi$.
When $\varepsilon<\theta<\pi-\varepsilon$ we have the uniform estimate

$$
\left|  \widehat{\mu_\beta}(\rho\omega^{\prime}\sin\theta,\rho\cos\theta)\right|
\lesssim\rho^{-\frac{n-1}{2}}.
$$
Hence, when $p>\frac{2(\beta-1)}{\beta-2},$

$$
\eqalign{
\int_{0}^{\pi}&\int_{\Bbb S^{n-2}}\left|  \widehat{\mu_\beta}(\rho\omega^{\prime
}\sin\theta,\rho\cos\theta))\right|  ^{p}\sin^{n-2}\theta
d\omega^{\prime}d\theta\cr 
&  \lesssim\,\rho^{-\frac{n-1}{2}p}\cr
&  +\int_0^{\varepsilon}\left(  \min\left(  \rho^{-\frac{n-1}{\beta
}},\rho^{-\frac{n-1}{2}}\left(  \sin\theta\right)  ^{-\frac{(n-1)(\beta
-2)}{2(\beta-1)}}\left(  \cos\theta\right)  ^{-\frac{(n-1)}{2(\beta-1)}
}\right)  \right)  ^{p}\left(  \sin\theta\right)  ^{n-2}d\theta\cr
&  \lesssim\,\rho^{-\frac{n-1}{2}p}+\int_0^{\varepsilon}\left(
\min\left(  \rho^{-\frac{n-1}{\beta}},\rho^{-\frac{n-1}{2}}\theta
^{-\frac{(n-1)(\beta-2)}{2(\beta-1)}}\right)  \right)  ^{p}\theta
^{n-2}d\theta\cr
&  \lesssim\rho^{-\frac{n-1}{2}p}+\int_0^{\rho^{-1+1/\beta}}
\rho^{-\frac{n-1}{\beta}p}\theta^{n-2}d\theta\cr
&  +\int_{\rho^{-1+1/\beta}}^{\varepsilon}\rho^{-\frac{n-1}{2}
p}\theta^{-\frac{(n-1)(\beta-2)p}{2(\beta-1)}}\theta^{n-2}d\theta\cr
&  \lesssim\,\rho^{-\frac{n-1}{2}p}+\rho^{-(n-1)\frac{p+\beta-1}{\beta}
}+\rho^{-(n-1)\frac{p+\beta-1}{\beta}}\lesssim\rho^{-(n-1)\frac{p+\beta
-1}{\beta}}}
$$

The computations when $p=\frac{2(\beta-1)}{\beta-2}$ or $p<\frac
{2(\beta-1)}{\beta-2}$ are similar.\qed
\medskip

Incorporating the average Fourier transform decay estimate of Prop. 2 into the
proof of [RT] as described in
\S2 then  yields Thm. 6.

\bigskip
\head{\bf \S4. Translates of transformations of a fixed surface}\endhead

We now turn to the proof of Thm. 4. If the number of parameters $m$ is greater than
$n-1$, then under the rank assumption of Thm. 4, we may select $n-1$ variables
$s_{i_1},\dots,s_{i_{n-1}}$ such that the corresponding square submatrix of
(1.4) is nonsingular. The estimate (1.5) then holds,  with
respect to $ds_{i_1}\dots ds_{i_{n-1}}$, uniformly in the other $s$ variables.
Since
$K\subset\subset\R^m$, we may integrate in all the variables and see that (1.5)
holds. Hence it sufffices to assume that $m=n-1$ and (1.4) is
nonsingular.

Starting with the two-dimensional case, let
$$
T_{s}:\R^{2}\rightarrow\R^{2}
$$
be a smooth family of transformations of the plane with $s\in\lbrack
a,b]$. Let $\gamma^0:[-\varepsilon,\varepsilon]\rightarrow\R^{2}$ be a
convex curve in $\R^{2}$ and let $\gamma_{s}\left(  t\right)
=T_{s}\left(  \gamma^0\left(  t\right)  \right)  $. We are interested in
the operator
$$
Rf\left(  x,s\right)  =\int f\left(  x+\gamma_{s}\left(  t\right)
\right)  dt=(f\ast\mu_{s})\left(  x\right)
$$
where $\mu_{s}$ is the measure defined by
$$
\int f\left(  x\right)  d\mu_{s}=\int f\left(  T_{s}(\g^0(t)) \right)  dt.
$$

Splitting the curve into a finite number of segments if necessary,  we may assume
the existence of two constants
$\varphi_{1}$ and $\varphi_{2}$ such that for any $t$ the left and the right
tangent lines at $\gamma(t)$ have slopes between $\varphi_{1}$ and
$\varphi_{2}$, and $\varphi_{2}-\varphi_{1}$ is small.

\bigskip
\proclaim{Proposition 3}
Let $\Phi^t=\left(  \cos\varphi,\sin\varphi\right)  $ with $\varphi_{1}
\leq\varphi\leq\varphi_{2}$. We assume that for every $\varphi_{1}\leq
\varphi\leq\varphi_{2}$ the matrix
$$
C=\Bigl[
\matrix
J_{T_{s}}\Phi^{t} \cr J_{\frac{dT_{s}}{ds}}\Phi^{t}
\endmatrix
\Bigr] \tag4.1
$$
is non-singular. Then
$$
\Bigl(\int\left\vert \widehat{\mu_{s}}\left(  \xi\right)  \right\vert
^{2}ds\Bigr)^\frac12\lesssim\left\vert \xi\right\vert ^{-\frac12}.\tag4.2
$$

\endproclaim
\bigskip

\fp{\bf Proof.} For any $\xi\neq0$ let
$$
\delta=\max_{s}\left(  \frac{\left\vert \xi\cdot J_{T_{s}}
\Phi\right\vert }{\left\vert \xi\right\vert },\frac{\left\vert \xi\cdot
J_{_{\frac{dT_{s}}{d\theta}}}\Phi\right\vert }{\left\vert \xi\right\vert
}\right).
$$
The nonsingularity of $C$ ensures that $\delta>c>0$. Let $\eta_{1}\left(
s\right)  ,\eta_{2}\left(  s\right)  \in C^{\infty}$ be cut-off
functions such that $\eta_{1}\left(  s\right)  +\eta_{2}\left(
s\right)  \equiv1$, $\frac{\left\vert \xi\cdot J_{\frac{dT_{s}
}{ds}}\Phi\right\vert }{\left\vert \xi\right\vert }>c$ on the support of
$\eta_{1}$ and $\frac{\left\vert \xi\cdot J_{_{T_{s}}}\Phi\right\vert
}{\left\vert \xi\right\vert }>c$ on the support of $\eta_{2}$. Then
$$
\int\left\vert \widehat{\mu_{s}}\left(  \xi\right)  \right\vert
^{2}ds=\int\left\vert \widehat{\mu_{s}}\left(  \xi\right)
\right\vert ^{2}\eta_{1}\left(  s\right)  ds+\int\left\vert
\widehat{\mu_{s}}\left(  \xi\right)  \right\vert ^{2}\eta_{2}\left(
s\right)  ds=I+II
$$

To estimate $I$ observe that
$$
I=\int\left\vert \widehat{\mu_{s}}\left(  \xi\right)  \right\vert
^{2}\eta_{1}\left(  s\right)  ds=\int\int\int e^{-2\pi i\xi
\cdot\left(  \gamma_{s}\left(  t\right)  -\gamma_{s}\left(
t'\right)  \right)  }\eta_{1}\left(  s\right)  ds dt' dt.
$$
Since
$$
\eqalign{
&  \left\vert \frac{d}{ds}\xi\cdot\left(  \gamma_{s}\left(
t\right)  -\gamma_{s}\left(  t'\right)  \right)  \right\vert =\left\vert
\xi\cdot\frac{dT_{s}}{ds}\left(  \gamma\left(  t\right)  \right)
-\xi\cdot\frac{dT_{s}}{ds}\left(  \gamma\left(  t'\right)  \right)
\right\vert \cr
&  =\left\vert \xi\cdot J_{\frac{dT_{s}}{ds}}\left(  c\right)
\left(  \gamma\left(  t\right)  -\gamma\left(  t'\right)  \right)  \right\vert
\geqslant c\left\vert \xi\right\vert \left\vert \gamma\left(  t\right)
-\gamma\left(  t'\right)  \right\vert}$$
(observe that the direction of $\gamma\left(  t\right)  -\gamma\left(
t'\right)  $ is between $\varphi_{1}$ and $\varphi_{2}$ by the convexity of
$\gamma)$), integrating by parts in $s$ yields
$$
I\lesssim\int\int\frac{1}{1+\left\vert \xi\right\vert \left\vert
\gamma\left(  t\right)  -\gamma\left(  t'\right)  \right\vert }dt'dt\lesssim
\frac{1}{\left\vert \xi\right\vert }.
$$
We now consider $II$. When $s$ belongs to the support of $\eta_{2}$ we
have, a.e. in $t$,
$$
\left\vert \frac{d}{dt}\xi\cdot T_{s}\left(  \gamma\left(  t\right)
\right)  \right\vert =\left\vert \xi\cdot J_{T_{s}}\left(  \gamma\left(
t\right)  \right)  \gamma^{\prime}\left(  t\right)  \right\vert \geqslant
c\left\vert \xi\right\vert
$$
and therefore, integrating by parts, we get
$$
\left\vert \widehat{\mu_{s}}\left(  \xi\right)  \right\vert =\left\vert
\int e^{-2\pi i\xi\cdot T_{s}\left(  \gamma\left(  t\right)  \right)
}dt\right\vert \lesssim\frac{1}{\left\vert \xi\right\vert },
$$
finishing the proof of Prop. 3.\qed
\bigskip

It is also possible to apply the geometric combinatorics technique of Christ[C]
(see also [TaW]) to obtain all but the sharp $L^{\frac32}\lra L^3$ result, with a
restricted weak-type estimate
at the endpoints, under the same geometric condition (4.1). 
Let $\gamma_{t,s}(x)= x-T_s(\gamma^0(t))$, thought of as a family of
diffeomorphisms of $\R^2$, indexed by $t,s$. Define, for
$y_0\in\R^2$ fixed, a map
$\Psi:\R^3\lra\R^2$ by
$$\Psi(t',t,s)=\gamma_{t',s}(\gamma^{-1}_{t,s}(y_0)).$$
Then the crucial things one needs for the argument of [C] to work are:
$$||D\Psi||_{2\times 2}\ge c|t'-t|,\tag{i}$$
where $||A||_{2\times 2}$ is the maximal $2\times 2$ minor of a $2\times
3$ matrix $A$, and
$$|\Psi^{-1}(x)|_1\le C,\tag{ii}$$
an upper bound on the lengths of the preimages of points under $\Psi$.
In this translation-invariant situation,
$$\Psi(t',t,s)= y_0+T_s(\gamma(t))-T_s(\gamma(t'))$$ 
so that
$$\eqalign{D\Psi =& \Bigl[ -DT_s(\dot\gamma(t')),\quad DT_s
(\dot\gamma(t)),\quad T'_s(\gamma(t))-T'_s(\gamma(t'))\Bigr]\cr
\simeq&\Bigl[-DT_s(\dot\gamma(s)),\quad
(t'-t)DT_s(\ddot\gamma(t'),\quad
(t'-t)DT'_s(\dot\gamma(t'))\Bigr],}$$
from which one sees that (i) will follow if
$$rank\Bigl[J_{T_s}\dot\gamma,\quad J_{T_s}\ddot\gamma,\quad
J_{\frac{dT_s}{ds}}\dot\gamma\Bigr]=2.$$ 
Since this approach does not yield the endpoint result, but only restricted weak
type, we shall not describe it in more detail.
\bigskip
In  dimensions $n\ge 3$, we need to impose a regularity condition on the surface in
addition to convexity, so we now assume that $\g$ is a $C^2$
convex parametrized hypersurface in $\R^n$.
By a partition of unity on the surface, we may assume that $\gamma\left(  t\right) 
$ is in a small neighborhood of a fixed point $x_{o}=\gamma\left(  t_{o}\right)  $.
We can also assume that the image of the Gauss map is a small compact subset
$\Omega$ of $\SN$, and  denote by $\Omega^{\perp}$ the set of directions that are
orthogonal to a direction in $\Omega$. 

\bigskip
\proclaim{Proposition 4}
Under the assumption that the matrix in (1.4) has rank $n$, 
$$
\left\{  \int\left\vert \widehat{\mu_{s}}\left(  \xi\right)  \right\vert
^{2}\Psi\left(  s\right)  ds\right\}  ^{\frac12}\lesssim\left\vert
\xi\right\vert ^{-\frac{n-1}{2}}\tag4.3
$$
where $\Psi\in C_{0}^{\infty}\left(  K\right)  $.
\endproclaim
\bigskip

\fp{\bf Remark.}
Let $\omega=\Xi^{t}J_{T_{s}}\left(  x_{0}\right)  $. By the rank assumption in
Thm. 4, we have that the vectors $\omega,\frac{\partial\omega}{\partial s_{1}
},\ldots,\frac{\partial\omega}{\partial s_{n-1}}$are linearly independent.
Hence, writing $s=(s',s'')\in\R^{n-1}\times\R^{m-n+1}$, it follows that for all
$s_0''\in\R^{m-n+1}$, the map
$s'\rightarrow
\omega(s',s_0'')$  defines a hypersurface in
$\R^{n}$ whose tangent hyperplane does not contain the origin. 
\medskip

\fp{\bf Proof.}
Let $\xi=\rho\Xi$ where $\rho=\left\vert \xi\right\vert $ and let
$$
\delta\left(  s\right)  =\inf_{\buildrel{\Phi\in\Omega^{\perp
}}\over{\left\vert \Phi\right\vert =1}}\max_{k=1,\ldots,n-1}\left\vert \Xi
^{t}J_{\frac{\partial T}{\partial s_{k}}}\left(  x_{o}\right)
\Phi\right\vert . \tag4.4
$$
By a smooth partition of unity we can assume that $\Psi$ is supported in a
small neighborhood of a fixed point $s_{0}$ and that on the support of
$\Psi$ either $\delta\left(  s\right)  \geq\varepsilon$
or
$\delta\left(  s\right)  <\varepsilon$
holds for a suitable $\varepsilon$ to be chosen later.

Assume we have $\delta\left(  s\right)  \geq\varepsilon$. Let us consider
$$
I=\int\left\vert \widehat{\mu_{s}}\left(  \xi\right)  \right\vert
^{2}\Psi\left(  t'\right)  ds=\int\int\int e^{-i\xi\cdot\left(
\gamma_{s}\left(  t\right)  -\gamma_{s}\left(  t'\right)  \right)
}\Psi\left(  t'\right)  ds dt' dt.
$$
We have
$$\eqalign{
&  \left\vert \frac{\partial}{\partial s_{k}}\xi\cdot\left(
\gamma_{s}\left(  t\right)  -\gamma_{s}\left(  t'\right)  \right)
\right\vert =\left\vert \xi\cdot\frac{\partial T_{s}}{\partial s_{k}
}\left(  \gamma\left(  t\right)  \right)  -\xi\cdot\frac{\partial T_{s}
}{\partial s_{k}}\left(  \gamma\left(  t'\right)  \right)  \right\vert \cr
&  =\rho\left\vert \Xi^t J_{\frac{\partial T_{s}}{\partial s_{k}}
}\left(  \gamma\left(  c\right)  \right)  \gamma^{\prime}\left(  c\right)
\left(  t-t'\right)  \right\vert \geqslant\varepsilon \rho\left\vert
\gamma^{\prime}\left(  c\right)  \left(  t-t'\right)  \right\vert \geq
c\varepsilon \rho\left\vert t-t'\right\vert .}
$$
We used the fact that $\gamma^{\prime}\left(  c\right)  \left(  t-t'\right)  $
is in the tangent hyperplane at $\gamma\left(  c\right)  $ and therefore is in
$\Omega^{\perp}$, and that $\gamma^{\prime}$ has maximal rank. Integrating by
parts $n-1$ times in $s_{k}$ we get
$$
I\lesssim\int\int\frac{1}{\left[  1+\rho\left\vert t-t'\right\vert \right]
^{n-1}}dsdt\lesssim \frac{1}{\rho^{n-1}}.
$$
We now consider the case
$$
\delta\left(  s\right)  <\varepsilon.\tag4.5
$$

Since the map $T_{s}\left(  x\right)  $ is smooth we can write
$T_{s}\left(  x\right)  =T_{s}\left(  x_{0}\right)  +J_{T_{s}
}\left(  x-x_{0}\right)  +E_{s}\left(  x,x_{0}\right)  $. Then
$$
\eqalign{
\left\vert \widehat{\mu_{s}}\left(  \xi\right)  \right\vert  &
=\left\vert \int e^{i\rho\Xi\cdot T_{s}\left(  \gamma\left(  t\right)
\right)  }dt\right\vert \cr
&  =\left\vert \int e^{i\rho\Xi\cdot\left[  J_{T_{s}}\left(  \gamma\left(
t\right)  -\gamma\left(  t_{0}\right)  \right)  +E_{s}\left(
\gamma\left(  t\right)  ,\gamma\left(  t_{o}\right)  \right)  \right]
}dt\right\vert \cr
&  =\left\vert \int e^{i\rho\Xi\cdot J_{T_{s}}\gamma\left(  t\right)
+i\rho\Xi\cdot E_{s}\left(  \gamma\left(  t\right)  ,\gamma\left(
t_{o}\right)  \right)  }dt\right\vert .}
$$
Setting $\omega=J_{T_{s}}^{t}\Xi$, 
one has
$$
\left\vert \widehat{\mu_{s}}\left(  \xi\right)  \right\vert =\left\vert
\int e^{i\rho\omega\cdot\gamma\left(  t\right)  +i\rho\Xi\cdot E_{s}\left(
\gamma\left(  t\right)  ,\gamma\left(  t_{o}\right)  \right)  }dt\right\vert
\tag4.6
$$

By (4.4) and (4.5) there exists $\Phi\in\Omega^{\perp}$ so
that for every $k$,
$$
\left\vert \Xi^{t}J_{\frac{\partial T}{\partial s_{k}}}\left(
x_{o}\right)  \Phi\right\vert <\varepsilon, \tag4.7
$$
i.e.,
$$
\left\vert \frac{\partial\omega}{\partial s_{k}}\cdot\Phi\right\vert
<\varepsilon. \tag4.8
$$

Without loss of generality we can assume that $\gamma\left(  t\right)
=\left(  t,\Gamma\left(  t\right)  \right)  $ with $\gamma^{\prime}\left(
t_{0}\right)  =0$ and that $\Phi=\left(  1,0,\ldots,0\right)  ^{t}$. We claim
that the Jacobian associated to the change of variables
$$
\eqalign{
\omega_{2}  &  =\omega_{2}\left(  s\right)\cr
&  \vdots\cr
\omega_{N}  &  =\omega_{N}\left(  s\right)}\tag4.9
$$
is nonsingular. Indeed, let $\omega^{\prime}=\left(  \omega_{2},\ldots
,\omega_{N}\right)  ^{t}$ and assume that the vectors
$$
\frac{\partial\omega^{\prime}}{\partial s_{1}},\ldots,\frac{\partial
\omega^{\prime}}{\partial s_{n-1}}
$$
are linearly dependent. Then for suitable $\left(  c_{1},\ldots
,c_{n-1}\right)  \neq0$ we have
$$
c_{1}\frac{\partial\omega^{\prime}}{\partial s_{1}}+\cdots+c_{n-1}
\frac{\partial\omega^{\prime}}{\partial s_{n-1}}=\left(  0,\ldots0\right)
^{t}.
$$
Let
$$
c_{1}\frac{\partial\omega_{1}}{\partial s_{1}}+\cdots+c_{n-1}
\frac{\partial\omega_{1}}{\partial s_{n-1}}=\alpha;
$$
then
$$
c_{1}\frac{\partial\omega}{\partial s_{1}}+\cdots+c_{n-1}\frac
{\partial\omega}{\partial s_{n-1}}=(\alpha,0,\ldots,0)^{t}.
$$
Since the vectors $\frac{\partial\omega}{\partial s_{1}},\ldots
,\frac{\partial\omega}{\partial s_{n-1}}$ are linearly independent by the
assumption 
(1.4), we have $\alpha\neq0$ and therefore $\Phi=\left(
1,0,\ldots0\right)  ^{t}$ is linearly dependent of $\frac{\partial\omega
}{\partial s_{1}},\ldots,\frac{\partial\omega}{\partial s_{n-1}}$.
This contradicts (4.8).

Also observe that, by (4.7) and (1.4),
$$
\omega_{1}=\omega\cdot\Phi=\Xi^{t}J_{T_{s}}\Phi
$$
is bounded away from zero.

Let us consider the integral in (4.6). In order to integrate by parts in
the $t_{1}$ variable we observe that
$$
\eqalign{
&  \frac{\partial}{\partial t_{1}}\left(  \left[  \omega^{\prime}t+\omega
_{n}\Gamma\left(  t\right)  \right]  +\Xi\cdot E_{s}\left(  \gamma\left(
t\right)  ,\gamma\left(  t_{o}\right)  \right)  \right) \cr
&  =\omega_{1}+\omega_{n}\frac{\partial\Gamma}{\partial t_{1}}+\frac{\partial
}{\partial t_{1}}\Xi\cdot E_{s}\left(  \gamma\left(  t\right)
,\gamma\left(  t_{o}\right)  \right)  .
}
$$
Since $\gamma^{\prime}\left(  t_{0}\right)  =0$, taking $t$ in a sufficiently
small neighborhood of $t_{0}$ we get that $\omega_{n}\frac{\partial\Gamma}{\partial
t_{1}}$ is small. Moreover the same can be shown for the last term since $\gamma$
is $C^{1}$ and $T_{s}$ is smooth. This ensures that the above derivative
is bounded away from zero. Hence,
$$
\eqalign{
&  \int e^{i\rho\left[  \omega^{\prime}\cdot t+\omega_{n}\Gamma\left(  t\right)
\right]  +i\rho\Xi\cdot E_{s}\left(  \gamma\left(  t\right)  ,\gamma\left(
t_{o}\right)  \right)  }dt\cr
&  =\frac{1}{\rho}\int e^{i\rho\left[  \omega^{\prime}\cdot
t+\omega_{n}\Gamma\left( t\right)  \right]  +i\rho\Xi\cdot E_{s}\left( 
\gamma\left(  t\right) ,\gamma\left(  t_{o}\right)  \right)  }H\left( 
\omega^{\prime},t\right)  dt,}
$$
where $H$ is a bounded function smooth in the $\omega^{\prime}$ variable. It
follows that
$$
\eqalign{
I^{2}  &  =\int\left\vert \widehat{\mu_{s}}\left(  \xi\right)
\right\vert ^{2}\Psi\left(  s\right)  ds\cr
&  =\frac{1}{\rho^{2}}\int\left\vert \int e^{i\rho\left[  \omega^{\prime}\cdot
t+\omega_{n}\Gamma\left(  t\right)  \right]  +i\rho\Xi\cdot E_{s}\left(
\gamma\left(  t\right)  ,\gamma\left(  t_{o}\right)  \right)  }H\left(
\omega^{\prime},t\right)  dt\right\vert ^{2}\Psi\left(  s\right)  ds
}
$$
Performing the change of variables (4.9), we obtain
$$
I^{2}=\frac{1}{\rho^{2}}\int\left\vert \int e^{i\rho\left[  \omega^{\prime}\cdot
t+\omega_{n}\Gamma\left(  t\right)  \right]  +i\rho\Xi\cdot E_{s}\left(
\gamma\left(  t\right)  ,\gamma\left(  t_{o}\right)  \right)  }H\left(
\omega^{\prime},t\right)  dt\right\vert ^{2}\widetilde{\Psi}\left(
\omega^{\prime}\right)  d\omega^{\prime}
$$
and thus

$$
I^{2}=\frac{1}{\rho^{2}}\int\left\vert \int\int e^{^{iR\left[  \omega^{\prime
}\cdot t+\omega_{n}\Gamma\left(  t\right)  \right]  +iR\Xi\cdot E_{s}\left(
\gamma\left(  t\right)  ,\gamma\left(  t_{o}\right)  \right)  }}H\left(
\omega^{\prime},t\right)  dt^{\prime}dt_{1}\right\vert ^{2}\widetilde{\Psi
}\left(  \omega^{\prime}\right)  d\omega^{\prime}%
$$
where $t^{\prime}=\left(  t_{2},\ldots,t_{n-1}\right) $. By the Minkowski
integral inequality we can bound $I$ by
$$
\frac{1}{\rho}\int\left\{  \int\left\vert \int e^{^{i\rho\left[  \omega
^{\prime\prime}\cdot t^{\prime}+\omega_{n}\Gamma\left(  t\right)  \right]
+i\rho\Xi\cdot E_{s}\left(  \gamma\left(  t\right)  ,\gamma\left(
t_{o}\right)  \right)  }}H\left(  \omega^{\prime},t\right)  dt^{\prime
}\right\vert ^{2}\widetilde{\Psi}\left(  \omega^{\prime}\right)
d\omega^{\prime}\right\}  ^{1/2}dt_{1}.
$$

Expanding and rewriting the term inside the brackets, we turn it into
$$
\eqalign{
&  \int\int\int e^{i\rho\left[  \omega^{\prime\prime}\cdot \left(  t^{\prime}
-u^{\prime}\right)  +\omega_{n}\left(  \Gamma\left(  t_{1},t^{\prime}\right)
-\Gamma\left(  t_{1},u^{\prime}\right)  \right)  \right]  +i\rho\Xi\cdot
E_{s}\left(  \gamma\left(  t\right)  ,\gamma\left(  t_{o}\right)
\right)  -i\rho\Xi\cdot E_{s}\left(  \gamma\left(  t_{1},u^{\prime}\right)
,\gamma\left(  t_{o}\right)  \right)  }\times\cr
&  \quad\quad\quad\quad H\left(  \omega^{\prime},t_{1},t^{\prime}\right)  H\left(
\omega^{\prime},t_{1},u^{\prime}\right)  \widetilde{\Psi}\left(
\omega^{\prime}\right)  d\omega^{\prime}dt^{\prime}du^{\prime}
}
$$

For fixed value of $t_{1},t^{\prime},u^{\prime}$, let $k$ be such that
$$
\left\vert t_{k}-u_{k}\right\vert >c\left\vert t-u\right\vert .
$$
Then the derivative of the phase in $\omega_{k}$ is controlled by
$$
\eqalign{
&  \left\vert t_{k}-u_{k}+\frac{\partial\omega_{n}}{\partial\omega_{k}}
\nabla\Gamma\cdot\left(  0,t^{\prime}-u^{\prime}\right)  +O\left(  t-u\right)
\max\left(  \left\vert t-t_{0}\right\vert ,\left\vert u-t_{0}\right\vert
\right)  \right\vert \cr
&  \geq\frac{1}{2}\left\vert t_{k}-u_{k}\right\vert .
}
$$
Therefore we can integrate by parts $n-2$ times in the above integral and we
get a term controlled by
$$
\int\int\frac{1}{\left(  1+\rho\left\vert t^{\prime}-u^{\prime}\right\vert
\right)  ^{n-2}}dt^{\prime}du^{\prime}\lesssim\rho^{-\left(  n-2\right)  }.
$$
Hence,
$$
I\lesssim\rho^{-\frac{n}{2}},
$$
which is better than we need.
\qed

\bigskip

Using Props. 3 and 4 in the proof of [RT] yields Thm. 4.

\bigskip
\head{\bf \S5. Nondegenerate generalized Radon transforms}\endhead

Let $X$ and $Y$ be smooth manifolds, with cotangent bundles $T^*X$ and $T^*Y$
having zero sections 0, and $C\subset (T^*X\backslash 0)\times (T^*Y\backslash 0)$
a canonical relation. Then $I^r(X,Y;C)$ is the class of Fourier integral operators
$F:\Cal E'(Y)\lra\Cal D'(X)$ of order $r\in\R$ associated with $C$. (We refer to
[H1,H2] for the background material on Fourier integral operators.) $L^2$ estimates
for Fourier integral operators associated with a canonical relation $C$ depend on
the structure of the projections $\pr:C\lra T^*Y$ and $\pl:C\lra T^*X$. Assume
$dim X= dim Y + m$. The optimal
$L^2$ estimates, with an operator $F\in I^{r-\frac{m}4}(X,Y;C)$
mapping
$L^2_{\alpha,comp}(Y)\lra L^2_{\alpha-r,loc}(X)$, hold under the assumption that
$\pr$ is a submersion (which guarantees that $\pl$ is an immersion), together with
the weak additional requirement that the spatial projections $\pi_X:C\lra X$ and
$\pi_Y:C\lra Y$ are submersions [H1;H2,Thm.25.3.8]; such canonical relations
$C$ are called {\it nondegenerate}. If we consider the special case of a
generalized Radon transform $R$ given by (1.6), we have $Y=\R^n$ and
$X=\R^{n+m}$. One can embed $R$ in an analytic family of operators by inserting
the factor $|\theta|^{-z}$ into the oscillatory representation (1.7); then
$$
R^z\in I^{-Re(z)-\frac{k}2-\frac{m}4}(\R^{n+m},\R^n;C)
$$
with $C=N^*Z'$, where $Z\subset\R^{n+m}\times\R^n$ is the incidence relation for
$R$, as explained in the Introduction. Under the assumption that $C$ is
nondegenerate, we have
$$
R^z: L^2(\R^n)\lra L^2(\R^{n+m}),\quad Re(z)=-\frac{k}2.\tag5.1
$$
On the other hand, the Schwartz kernel of $R^z$ is in $L^\infty$ for $Re(z)=n-k$,
so that
$$
R^z: L^1(\R^n)\lra L^\infty(\R^{n+m}),\quad Re(z)=n-k.\tag5.2
$$
The bounds in both (5.1) and (5.2) grow at most exponentially in $|Im(z)|$ and
hence Thm. 5 follows by analytic interpolation.

Next, we make the connection between Thms. 4 and 5 by showing that, if the
hypersurface $\gamma$ in Thm. 4 is $C^\infty$, then condition (1.4) implies that
the associated canonical relation $C$ is nondegenerate, so that in this case Thm.
4 becomes a special case of Thm. 5.  Note that, by a simple calculation,  (1.4)
holds for the family
$\{T_s\}$ iff it holds for the family $\{T_s^{-1}\}$, and for convenience we will
work with the latter. The support of the Schwartz kernel of $R$ is then
$$Z=\{(x,s,y): x-y\in T_s^{-1}(\gamma)\}.$$
Letting $\nu(t)$ be a unit normal at $t\in\gamma$, we have

$$\eqalign{C=&N^*Z'\quad \cr
=&\{(x,s,\theta\cdot J_{T_s}^*\nu(t),\theta\cdot
(\frac{\partial T_s}{\partial s})^*\nu(t); x-t,\theta\cdot J_{T_s}^*\nu(t))\cr
&\quad :
x\in\R^n,s\in
\R^m,t\in\gamma,\theta\in\R\backslash 0\}\cr
=&\{(*,*; y,\theta\cdot J_{T_s}^*\nu(t))
:y\in\R^n,s\in
\R^m,t\in\gamma,\theta\in\R\backslash 0\},}$$
from which we see that $rank(D\pi_R)=n+rank(\frac{D(\eta)}{D(s,t,\theta)})$. 
Condition (1.4) then implies that
$rank(\frac{D(\eta)}{D(\theta,s_1,\dots,s_{n-1})})=n$.

We also point out that Seeger[S] has obtained $L^p-L^q$ estimates for  generalized
Radon transforms, almost sharp in the finite-type setting in two dimensions. If
$C$ is nondegenerate,one can see that the 
$Z$ is of type (1,1) in the terminology of [S], and thus $R:L_{comp}^p(Y)\lra
L_{loc}^q(X)$ for 
$$(\frac1p,\frac1q)\in int\{hull((0,0),(1,1),(\frac23,\frac13))\}.$$
However, the commutator approach of [S] is insensitive to the presence of more
variables and therefore, in the particular context of Thm. 5, does not yield
estimates outside of this set, regardless of the dimension.

It is possible to formulate geometric criteria under which the canonical relation
$C$ is nondegenerate. First consider the case of curves, $k=1$. Write
$\g_t(\cdot,s)=\g(\cdot,s;t)$, so that $\{\g_t\}$ is a one-parameter family of
diffeomorphisms of $\R^n$ (parametrized by
$s\in\R^m$); by a change of variables in $\R^n$, we may assume $\g_0=Id$. As
described in [GS,Eqn. 6.5] (see also [CNSW,\S9.3]), we can parametrize
$C$ as 
$$C=\{(\g_t^{-1}(y;s),s,(D_x\g)^*(\eta),(D_s\g)^*(\eta);y,\eta):
(s,y,t)\in\R^{m+n+1},\eta\perp \G(y,s;t)\}$$
where $\G$ is the (right) pullback of $\dot\g$ by the family of
diffeomorphisms $\{\g_t\}$, namely
$$\Gamma(y,s;t)=\frac{d}{dt'}\bigl(\gamma_{t'+t}\circ\gamma_t^{-1}(y)\bigr)
\bigr|_{t'=0}.\tag5.3$$
For each $y,s$, $\G(y,s;\cdot):\R\lra T_y\R^n$. 
\bigskip
\fp {\bf Ex. 1}. If $\g(x,s;t)=x+\g^0(s;t)$ is a translation-invariant
family, then $\G(y,s;t)=\dot\g^0(s;t)$ is just the velocity vector of the
curve at time $t$.
\medskip
\fp {\bf Ex. 2}. If, as in [CNSW,\S9.1], we prescribe a variable family of
curves via a Taylor expansion in $t$ and the exponential map (and allow
$s$-dependence), 
$$\g(x,s;t)=exp_x(tX(s)+t^2Y(s)+\dots)$$
where $X(s), Y(s),\dots$ are vector fields on $\R^n$ depending on
$s\in\R^k$, then, as calculated in [GS,\S6.4], 
$$\G(y,s;t)=X(s)+2tY(s)+\dots,\tag5.4$$
which is enough to determine whether $C$ is nondegenerate.
\medskip

If we work locally in $x,s$ and $t$, so that the first component
$\dot\g_1\ne 0$, then $\G_1\ne 0$ as well (for $|t|$ small). Writing
$\eta=(\eta_1,\eta')$, etc., we may then solve $\eta\perp\G(y,s;t)$ for
$\eta_1$ in terms of $\eta'\in\R^{n-1}\back 0$ and write the projection
$\pr:C\lra \tyn$ as
$$\pr(s,y,t,\eta')=(y,-\frac{\G'(y,s;t)\cdot\eta'}{\G_1(y,s;t)},\eta').\tag5.5$$
The canonical relation $C$ is nondegenerate if $\pr:C\lra\tyn$ is a
submersion, which implies also that $\pl:C\lra\txn$ is an immersion. 
Thus, we have

\proclaim{Theorem 7} If $\g(x,s;t)$ is a $C^\infty$ family of curves in
$\R^n$ such that
$$\R^{k+1}\owns (s,t)\lra \frac{\G'(y,s;t)\cdot\eta'}{\G_1(y,s;t)}\hbox{
has no critical points }\forall\eta'\in\R^{n-1}\back 0\tag5.6$$
then 
$$R:L^{\frac{2n-1}{n}}(\R^n)\lra
L^{\frac{2n-1}{n-1}}(\R^n\times\R^m).\tag5.7$$
\endproclaim

\bigskip
Condition (5.6) can be restated as a maximal rank condition on a $(m+1)\times
(n-1)$ matrix:
$$rank\Bigl[\matrix \G_1\cdot
D_s\G'-D_s\G_1\otimes\G'\cr \G_1\cdot
\partial_t\G'-\partial_t\G_1\otimes\G'\endmatrix\Bigr]=n-1.$$
Thus, a necessary condition for $C$ to be nondegenerate is that $m\ge
n-2$.

For the translation-invariant Ex. 1 above, we may write
$\g^0(s;t)=(t,g(s;t))$, where $g:\R^{m+1}\lra \R^{n-1}$; then
$\G(y,s;t)=\dot\g^0=(1,\dot g(s;t))$, so that (5.6) becomes
$$rank\Bigl[\matrix  D_s\dot g\cr \ddot g\endmatrix\Bigr]=n-1.\tag5.8$$
\bigskip
\fp For $n=2,m=0$ (i.e., no $s$ parameter) we need $\ddot g\ne 0$ as in
the result of Littman[L] and Strichartz[Str], while for $n=2,m=1$, (5.8) becomes:
$\ddot g\ne 0$ or $\partial_s\dot g\ne 0$, which includes the result of [RT] in the
smooth setting.
In $\R^3$, we need at least $m=1$, and then (5.8) becomes
$\ddot g\wedge\partial_s\dot{g}\ne 0$, i.e., $\{\dot g,\ddot g,\partial_s\dot
g\}$ linearly independent. If the family $\g^0(s;\cdot)$ arises from
rotation of an initial curve $\g_0(\cdot)$ about an axis $\R\cdot v$,
then we need $\dot\g_0\wedge v\ne 0, \ddot\g_0\cdot v\ne 0$. For example, 
convolution with the rotations of $\g_0(t)=(t,t^2,0)$ about the $x_2$ axis in
$\R^3$ already maps $L^{\frac53}(\R^3)\lra L^{\frac52}(\R^3\times S^1)$.

Thm. 3 for curves ($k=1$) follows from Thm. 7, since we may take
$s\in\R^{\frac{n(n-1)}2}$ to be local coordinates on $SO(n)$ and (5.8)
holds; essentially this says that $SO(n)$ acts transitively on the
sphere.
\bigskip
If one wants to formulate the results in terms of averages over
$m$-dimensional families of $k$-surfaces in $\R^n$, then only a few
changes are necessary. Starting with a $C^\infty$ map
$$\g:\R^n\times\R^m\times\R^k\lra \R^n,\quad \g(x,s;0)=x,\quad
D_t\g\hbox{ an injection,}$$
the resulting generalized Radon transform belongs to
$I^{-\frac{k}2-\frac{m}4}(\R^{n+m},\R^n;C)$. To describe the canonical relation
$C$, we use the pullback
$$\G(y,s;t)=D_{t'}(\g_{t+t'}\circ\g_t^{-1})|_{t'=0},$$
which is a map $\G:\R^n\times\R^m\times\R^k\lra \R^{k^*}\otimes T_y\R^n$.
We can assume that, with $x=(x',x'')\in\R^k\times\R^{n-k}$, etc., we have
that
$D_{x'}\g'$ is nonsingular, and thus $\G'$ is nonsingular for $|t|$
small. Condition (5.6) is then replaced by 
$$rank(D_{s,t}(({\G'}^*)^{-1}(\G'')(\eta'')))=k,\quad\forall\eta''\in\R^{n-k}\back
0.\tag5.9$$
Under this assumption, $C$ is nondegenerate.
Again, specializing to the translation-invariant case and letting
$s\in\R^{\frac{n(n-1)}{2}}$ be local coordinates on
$SO(n)$, it is not hard to see that $(5.9)$ is satisfied for any smooth
initial $k$-surface, and thus Thm. 3 follows.

\end